\documentclass[12pt, reqno]{amsart}
\usepackage{amssymb}
\usepackage{amsfonts}
\usepackage{amssymb,amsmath,amscd}
\usepackage[colorlinks=true]{hyperref}
\usepackage{mathrsfs}
\newtheorem{Theorem}{\indent Theorem}[section]

\newtheorem{Lemma}[Theorem]{\indent Lemma}

\theoremstyle{remark}

\begin{document}
\centerline{
\bf Averages of exponential twists of the  von Mangoldt function
}
\bigskip
\centerline{Xiumin Ren\footnote{School of Mathematics,
               Shandong University,
               Jinan, Shandong, 250100, China, xmren@sdu.edu.cn} and Wei Zhang\footnote{School of Mathematics and Statistics, Henan University
               Kaifeng,   Henan 475004
               China, zhangweimath@126.com}}

\bigskip

\textbf{Abstract}
In this paper, we obtain some improved results for the exponential sum $\sum_{x<n\leq 2x}\Lambda(n)e(\alpha k n^{\theta})$ with $\theta\in(0,5/12),$ where $\Lambda(n)$ is the von Mangoldt function.  Such  exponential sums have relations with the so-called quasi-Riemann hypothesis and were considered by Vinogradov \cite{Va} and Murty-Srinivas \cite{Mu}.
\medskip

\medskip


\bigskip
\bigskip
\numberwithin{equation}{section}

\section{Introduction}
In this paper,  we are interested in the exponential sum
\[
S(k,x,\theta):=\sum_{x<n\leq 2x}\Lambda(n)e(k\alpha n^{\theta}),
\]
where $x\geq2$ and $k\in\mathbb{Z}^{+}$ are the main parameters, $\alpha\neq 0$ and $0<\theta<1$
are fixed, $\Lambda(n)$ is the von Mangoldt function, and $e(z) = e^{2\pi iz}.$

We
call
$S(k, x,\theta)$ Vinogradov's exponential sum, since
 it was first
considered by I. M.
Vinogradov \cite{Va} in the special case $\theta = 1/2$. Actually, he proved in \cite{Va} that, for $k \leq x^{1/10},$
\[ S(k, x, 1/2) \ll k^{1/4}x^{7/8+\varepsilon},
\]
where the implied
constant may depend on $\alpha$ and $\varepsilon$. Iwaniec
 and Kowalski (see
(13.55) in  \cite{IK}) remarked that the stronger inequality
\[S(1,x,1/2) \ll x^{5/6}\log^{4} x
\]
follows from an application of Vaughan's identity. For general $\theta$ and $k$, Murty and Srinivas \cite{Mu}  proved  that
\[S(k, x,\theta) \ll k^{1/8}x
^{(7+\theta)/8}\log (xk^{3}),\]
where the implied
constant may depend on $\alpha$ and $\theta$.
In 2006, Ren \cite{Ren} proved that
\begin{align}\label{Ren1}
S(k, x,\theta)\ll \left(k^{1/2}x^{(1+\theta)/2}+
x^{4/5}+k^{-1/2}x^{1-\theta/2}\right)\log ^{A}x,
\end{align}
for arbitrary $A>0$, and for $\theta\leq1/2$ and $k<x^{1/2-\theta},$
\begin{align}\label{Ren}
S(k, x,\theta)\ll\left( k^{1/10}x^{3/4+\theta/10}
+k^{-1/2}x^{1-\theta/2}\right)\log ^{11}x.
\end{align}

In this paper, we will prove the following Theorem 1.1, which is new for $\theta\in(0,5/12).$ In \cite{Iw},
Iwaniec, Luo and Sarnak showed that such type exponential sums are connected to the qausi-Riemann Hypothesis (or the existence of zero-free region) for $L(s,f)$, where $f$ is any holomorphic cusp form of integral weight for $SL(2,\mathbb{Z})$.

\begin{Theorem}\label{th1}
For $0<\theta<5/12$ and $1\leq k<x^{5/12-\theta-\varepsilon},$ there exists an absolute constant $c_0>0$ such that
\[
S(k, x,\theta)\ll k^{-1/2}x^{1-\theta/2}\exp(-c_0(\log x)^{1/3-\varepsilon}),
\]
where the implied constant may depend on $\alpha,$  $\theta$, and $\varepsilon$,  which denotes an arbitrarily small positive constant.
\end{Theorem}
Obviously, when $\theta<5/12$ and $k<x^{5/12-\theta-\varepsilon},$ Theorem 1.1 improve (\ref{Ren}). Some much sharper estimates can be obtained if one assumes the zero-density hypothesis, i.e,
\begin{align}\label{ZD}
N(\sigma,T)\ll T^{2(1-\sigma)}\log^{B} T,\ \ \sigma\geq1/2,
\end{align}
where $N(\sigma, T)$ is the number of zeros of $\zeta(s)$ in the
region $\{\sigma\leq \Re s\leq 1, |t|\leq T\}$ and $B$ is some positive constant.
In fact, under (\ref{ZD}),  it is proved in \cite{Ren} that
\begin{align}\label{ZDH}
S(k, x,\theta)\ll \left(k^{1/2}x^{(1+\theta)/2}+
k^{-1/2}x^{1-\theta/2}\right)\log ^{B+2}x,
\end{align}
where the implied constant may depend on $\alpha$, $\theta$ and $B$.
In fact, our idea can also be used to give a better result than (\ref{ZDH}).
\begin{Theorem}\label{RH}
Under {\rm (\ref{ZD})}, for $0<\theta<1/2$ and $1\leq k<x^{1/2-\theta-\varepsilon},$ there exists an absolute constant $c_0$ such that
\[
S(k, x,\theta)\ll k^{-1/2}x^{1-\theta/2}\exp(-c_0(\log x)^{1/3-\varepsilon}),
\]
where the implied constant may depend on $\alpha,$ $\varepsilon$ and  $\theta$.
\end{Theorem}

It is worth pointing out that, comparing with Theorem \ref{th1},  the ranges of $\theta$ and $k$ have been extended in Theorem \ref{RH}.
\medskip

\section{Proof of Theorem \ref{th1}}

To prove Theorem \ref{th1}, we will borrow the idea in \cite{Ren}  by using the results related to zeros of Riemann zeta function.
The following lemma will be used in the proof of Theorem \ref{th1} and Theorem \ref{RH}.
\medskip

\begin{Lemma}[see page 71 of \cite{T}]
\label{le3.1}
Let $F(u)$ and $G(u)$ be real functions in $[a,b],$ satisfying
$|G(u)|\leq M$ and that $G(u)$ and $1/F'(u)$ are monotone.\\
{\rm (1)}  If $F'(u)\geq m>0$ or $F'(u)\leq -m<0,$ then
\[
\int_{a}^{b}G(u)e(F(u))du\ll \frac{M}{m};
\]
{\rm (2)} If $F''(u)\geq r>0$ or $F''(u)\leq -r<0,$ then
\[
\int_{a}^{b}G(u)e(F(u))du\ll \frac{M}{\sqrt{r}}.
\]
\end{Lemma}

{\bf Proof of Theorem \ref{th1}}
Using partial summation and the explicit formula (see  (5.53) in \cite{IK}): for $1\leq T\leq x$,
\[
\sum_{n\leq x}\Lambda(n)=x-\sum_{|\gamma|\leq T}\frac{x^{\rho}}{\rho}+O\left(\frac{x}{T}(\log xT)^{2} \right),
\]
we have
\begin{align}\label{30.1}
\sum_{x<n\leq2x}\Lambda(n)e( k\alpha n^{\theta})
&=\int_{x}^{2x}e(k\alpha u^{\theta})d\sum_{n\leq u}\Lambda(n)\nonumber\\
&=\int_{x}^{2x}e(k\alpha u^{\theta})du-\sum_{|\gamma|\leq T}\int_{x}^{2x}u^{\rho-1}e(k\alpha u^{\theta})du\nonumber\\
&+O\left((1+ k|\alpha| x^{\theta})\frac{x \log ^{2}x}{T} \right).
\end{align}
Here $\rho=\beta+i\gamma$ denotes a  zero of $\zeta(s)$ with  $0<\beta<1$, $|\gamma|\leq T$.
Set
\[
T=T_{0}=x,
\]
then the error-term is  $O(((1+k|\alpha|x^{\theta})\log^2x) =O_{\alpha}(kx^\theta\log^2x)$.
Moreover, we  have
\begin{align}\label{30.2}
\int_{x}^{2x}e(k\alpha u^{\theta})du=
\frac1\theta \int_{x^{\theta}}^{(2x)^{\theta}}u^{1/\theta-1}e(k\alpha u)du \ll_{\alpha,\theta} k^{-1}x^{1-\theta}.
\end{align}
Making the change of variable $u^{\theta}=v,$ we get
\begin{align*}
\int_{x}^{2x}u^{\rho-1}e(k\alpha u^{\theta})du=\frac1{\theta}\int_{x^{\theta}}^{(2x)^{\theta}}v^{\frac\beta\theta-1}e(f(v))dv,
\end{align*}
where
$$
f(v)= k\alpha v +\frac{\gamma}{2\pi\theta}\log v.
$$
Trivially one has
\begin{align}\label{33.1}
\int_{x}^{2x}u^{\rho-1}e(k\alpha u^{\theta})du
\ll x^\beta.
\end{align}
On the other hand we have
\begin{align*}
&|f'(v)|=|k\alpha+\frac{\gamma}{2\pi\theta v}|\geq \frac{\min_{v\in[x^{\theta},(2x)^{\theta}]}|\gamma
+2\theta\pi k\alpha v|}{2\pi\theta |v|},\\
&|f''(v)|=\frac{|\gamma|}{2\pi\theta v^{2}}.
\end{align*}
By Lemma \ref{le3.1} and (\ref{33.1}) we get
\begin{align*}
\int_{x^{\theta}}^{(2x)^{\theta}}v^{\frac\beta\theta-1}e\left(f(v)\right)dv\ll\begin{cases}
\frac{x^{\beta}}{\sqrt{1+\theta k|\alpha| x^{\theta}}}     &\textup{for}\ \  |\gamma| \leq 4(1+\theta\pi k|\alpha|(2x)^{\theta}),\\
\\
\frac{x^{\beta}}{ 1+|\gamma|}      &\textup{for}\ \  4(1+\theta\pi k|\alpha|(2x)^{\theta})\leq|\gamma| \leq T_{0}.
\end{cases}
\end{align*}
Therefore
\begin{align*}
&\sum_{|\gamma|\leq T}\int_{x}^{2x}u^{\rho-1}e(k\alpha u^{\theta})du\\
&\ll\frac{1}{\sqrt{1+\theta k|\alpha|x^{\theta}}}\sum_{|\gamma| \leq 4(1+\theta\pi k|\alpha|(2x)^{\theta})}x^{\beta}+
\sum_{4(1+\theta\pi k|\alpha|(2x)^{\theta})\leq|\gamma| \leq T_0}\frac{x^{\beta}}{1+|\gamma|}.
\end{align*}

%\begin{Lemma}Let $\varepsilon>0$ and $T=y^{\theta},$ where $\theta<5/12.$ Then
%\[
%\sum_{\gamma\leq T}y^{\beta}\ll y \exp\left(-c(\log y)^{1/3+\varepsilon}\right)
%\]
%where the implied constant depends on $\theta$ and $\varepsilon$.
%\end{Lemma}
%\begin{proof}
%Let $\delta= 5/12-\theta>0$ and $\sigma_{0}=1-c\left(\log T \right)^{-2/3}\left(\log \log T \right)^{-1/3}.$ Then
%\begin{align*}
%\sum_{|\gamma\leq T|}y^{\beta}&\ll (\log y)\sup_{1/2\leq \sigma\leq1}N(\sigma,y^{\theta})y^{\sigma}\\
%&\ll (\log y)^{C} \sup_{1/2\leq \sigma\leq\sigma_{0}}y^{\sigma+\theta(12/5)(1-\sigma)}\\
%&\ll y(\log y)^{C}\sup_{1/2\leq \sigma\leq\sigma_{0}}y^{-\delta(1-\sigma)}\\
%&\ll y(\log y)^{C}\sup_{1/2\leq \sigma\leq\sigma_{0}}y^{-c'(\log y)^{-2/3}(\log \log y)^{-1/3}}\\
%&\ll y\exp(-c(\log y)^{-1/3+\varepsilon}).
%\end{align*}
%\end{proof}
Assume that, for some positive constant $C$,
\begin{align*}
N(\sigma,T)\ll T^{A(\sigma)(1-\sigma)}\log^{C} T.
\end{align*}
Then by the Riemann-Von Mangoldt formula, for $2\leq U\leq T_0$ we have
\begin{align*}
\sum_{|\gamma|\leq U}x^{\beta}&=-\int^{1}_{0}x^{\sigma}dN(\sigma,U)\nonumber\\
&\ll x^{1/2}U\log U+(\log U)^{C}\log x\max_{1/2\leq\sigma\leq\sigma_0}U^{A(\sigma)(1-\sigma)}x^{\sigma},
\end{align*}
where
\[
\sigma_{0}=1-c_0\left(\log T \right)^{-2/3}\left(\log \log T \right)^{-1/3}
\]
with $c_0$ an absolute positive constant.
Here we have used the well known zero-free region results (for example, see \cite{IK,T}) which states that $\zeta(s)\not=0$ for
$
\sigma>\sigma_0,
$

Let $x$ be sufficiently large such that $\theta\pi k|\alpha| (2x)^{\theta}\gg 1$,  then we have
\begin{align*}
&\frac{1}{\sqrt{1+\theta k|\alpha| x^{\theta}}}\sum_{|\gamma| \leq 4(1+\theta\pi k|\alpha|(2x)^{\theta})}x^{\beta}\nonumber\\
&\ll (\log x)^{C+1}\big(k^{1/2}x^{(1+\theta)/2}
+\max_{1/2\leq\sigma\leq\sigma_0}k^{A(\sigma)(1-\sigma)
-1/2}x^{\sigma+\theta A(\sigma)(1-\sigma)-\theta/2}\big),
\end{align*}
and
\begin{align*}
&\sum_{4(1+\theta\pi k|\alpha|(2x)^{\theta})\leq|\gamma| \leq T_0}\frac{x^{\beta}}{1+|\gamma|}\\
&\ll (\log x)\max_{4(1+\theta\pi k|\alpha|(2x)^{\theta})<T_{1}\leq T}T_{1}^{-1}
\sum_{T_{1}\leq|\gamma|\leq2T_{1}}x^{\beta}\\
&\ll (\log x)^{C+2}\big(x^{1/2}+\max_{1/2\leq\sigma\leq\sigma_0}
k^{A(\sigma)(1-\sigma)-1}x^{\sigma+\theta A(\sigma)(1-\sigma)-\theta}\big).
\end{align*}
Writing
\[
g(\sigma)=\sigma+\theta A(\sigma)(1-\sigma)-\frac{\theta}{2},
\]
and collecting the above estimates we  get
\begin{align*}
&\sum_{|\gamma|\leq T}\int_{x}^{2x}u^{\rho-1}e(k\alpha u^{\theta})du\\
&\ll (\log x)^{C+2}\big( k^{1/2}x^{(1+\theta)/2}
+\max_{1/2\leq\sigma\leq\sigma_0}k^{A(\sigma)(1-\sigma)
-1/2}
x^{g(\sigma)}\big).
\end{align*}

By the well know result of Ingham \cite{In} and Huxley \cite{Hu}, we can choose $A(\sigma)=12/5.$
Thus we have
\begin{align*}
&\max_{1/2\leq\sigma\leq \sigma_{0}}
k^{A(\sigma)(1-\sigma)-1/2}x^{g(\sigma)}\\
&\ll (\log x)^{C_1} \sup_{1/2\leq \sigma\leq\sigma_{0}}k^{A(\sigma)(1-\sigma)-1/2}
x^{\sigma+12\theta (1-\sigma)/5-\theta/2}\\
&\ll k^{-1/2}x^{1-\theta/2}(\log x)^{C_1}\sup_{1/2\leq \sigma\leq\sigma_{0}}\left(k^{12/5}x^{12\theta/5-1}\right)^{1-\sigma}.
\end{align*}
Thus for $\theta<5/12$ and $k<x^{5/12-\theta-\varepsilon},$ we get
\begin{align*}
&\max_{1/2\leq\sigma\leq \sigma_{0}}
k^{A(\sigma)(1-\sigma)-1/2}x^{g(\sigma)}\\
&\ll k^{-1/2}x^{1-\theta/2}(\log x)^{C_1}\sup_{1/2\leq \sigma\leq\sigma_{0}}x^{-c_0(\log x)^{-2/3}(\log \log x)^{-1/3}}
\\
&\ll k^{-1/2}x^{1-\theta/2}\exp(-c_0(\log x)^{1/3} (\log x\log x)^{-1/3})\\
&\ll k^{-1/2}x^{1-\theta/2}\exp(-c_0(\log x)^{1/3-\varepsilon}).
\end{align*}
\medskip
This together with  (\ref{30.1}) and (\ref{30.2}) shows that, for $\theta\in(0,5/12)$ and $1\leq k<x^{5/12-\theta-\varepsilon}$,
\begin{align*}
&\sum_{x<n\leq2x}\Lambda(n)e( \alpha n^{\theta})\\
&\ll
k^{1/2}x^{(1+\theta)/2}(\log x)^{C}
+k^{-1/2}x^{1-\theta/2}\exp(-c_0(\log x)^{1/3-\varepsilon})+k^{-1}x^{1-\theta}+kx^{\theta}\\
&\ll k^{-1/2}x^{1-\theta/2}\exp(-c_{0}(\log x)^{1/3-\varepsilon}).
\end{align*}
This finishes the proof of Theorem 1.1. $\square$

\bigskip

$\mathbf{Acknowledgement}$ This work was supported by National Natural Science Foundation of China (Grant No. 11871307).

\end{document}